\documentclass[12pt,leqno]{article}
\usepackage{amsthm,amsbsy,amsfonts,%
   amssymb,amsmath,amscd}
\usepackage{latexsym}
\usepackage{fontenc}
\usepackage[mathscr]{euscript}

\newcommand{\Q}{\mathbb Q}
\newcommand{\Z}{\mathbb Z}
\newcommand{\R}{\mathbb R}
\newcommand{\C}{\mathbb C}

\newcommand{\F}{\mathbb F}

\newcommand{\A}{\mathbb A}

\newcommand{\gQ}{{\rm Gal}(\overline \Q/\Q)}
\newcommand{\scR}{\mathcal R}
\newcommand{\scA}{\mathcal A}
\newcommand{\scO}{\mathcal O}
\newcommand{\scG}{\mathcal G}

\title{Irreducibility and cuspidality}
\author{Dinakar Ramakrishnan\footnote{Partially supported
by the NSF through the grant DMS-0402044}}
\date{}

\begin{document}

\maketitle

\medskip


\bigskip

\section*{Introduction}

\bigskip

{\it Irreducible} representations are the building blocks of
general, semisimple Galois representations $\rho$, and {\it
cuspidal} representations are the building blocks of automorphic
forms $\pi$ of the general linear group. It is expected that when an
object of the former type is associated to one of the latter type,
usually in terms of an identity of $L$-functions, the irreducibility
of the former should imply the cuspidality of the latter, and
vice-versa. It is not a simple matter - {\it at all} - to prove this
expectation, and nothing much is known in dimensions $>2$. We will
start from the beginning and explain the problem below, and indicate
a result (in one direction) at the end of the Introduction, which
summarizes what one can do at this point. The remainder of the paper
will be devoted to showing how to deduce this result by a synthesis
of known theorems and some new ideas. We will be concerned here only
with the {\sl so called} {\it easier} direction of showing the
cuspidality of $\pi$ given the irreducibility of $\rho$, and refer
to [Ra5] for a more difficult result going the other way, which uses
crystalline representations as well as a refinement of certain deep
modularity results of Taylor, Skinner-Wiles, et al. Needless to say,
{\it easier} does not mean {\sl easy}, and the significance of the
problem stems from the fact that it does arise (in this direction)
naturally. For example, $\pi$ could be a functorial, automorphic
image $r(\eta)$, for $\eta$ a cuspidal automorphic representation of
a product of smaller general linear groups: $H(\A)=\prod_j
GL(m_j,\A)$, with an associated Galois representation $\sigma$ such
that $\rho=r(\sigma)$ is irreducible. If the automorphy of $\pi$ has
been established by using a flexible converse theorem ([CoPS1]),
then the cuspidality of $\pi$ is not automatic. In [RaS], we had to
deal with this question for cohomological forms $\pi$ on GL$(6)$,
with $H={\rm GL}(2)\times {\rm GL}(3)$ and $r$ the Kronecker
product, where $\pi$ ia automorphic by [KSh1]. Besides, the main
result (Theorem A below) of this paper implies, as a consequence,
the cuspidality of $\pi=$ sym$^4(\eta)$ for $\eta$ defined by any
non-CM holomorphic newform $\varphi$ of weight $\geq 2$ relative to
$\Gamma_0(N)\subset {\rm SL}(2,\Z)$, without appealing to the
criterion of [KSh2]; here the automorphy of $\pi$ is known by [K]
and the irreducibility of $\rho$ by [Ri].

Write $\overline \Q$ for the field of all algebraic numbers in $\C$,
which is an infinite, mysterious Galois extension of $\Q$. One could
say that the central problem in algebraic Number theory is to
understand this extension. {\it Class field theory}, one of the
towering achievements of the twentieth century, helps us understand
the {\it abelian} part of this extension, though there are still
some delicate, open problems even in that well traversed situation.

Let $\scG_\Q$ denote the absolute Galois group of $\Q$, meaning
$\gQ$. It is a profinite group, being the projective limit of finite
groups ${\rm Gal}(K/\Q)$, as $K$ runs over number fields which are
normal over $\Q$. For fixed $K$, the Tchebotarev density theorem
asserts that every conjugacy class $C$ in Gal$(K/\Q)$ is the {\it
Frobenius class} for an infinite number of primes $p$ which are
unramified in $K$. This shows the importance of studying the {\it
representations} of Galois groups, which are intimately tied up with
conjugacy classes. Clearly, every $\C$-representation, i.e., a
homomorphism into GL$(n,\C)$ for some $n$, of Gal$(K/\Q)$ pulls
back, via the canonical surjection $\scG_\Q \to {\rm Gal}(K/\Q)$, to
a representation of $\scG_\Q$, which is continuous for the profinite
topology. Conversely, one can show that every {\it continuous}
$\C$-representation $\rho$ of $\scG_\Q$ is such a pull-back, for a
suitable finite Galois extension $K/\Q$. E.~Artin associated an
$L$-function, denoted $L(s,\rho)$, to any such $\rho$, such that the
arrow $\rho\to L(s,\rho)$ is additive and inductive. He conjectured
that for any non-trivial, irreducible, continuous
$\C$-representation $\rho$ of $\scG_\Q$, $L(s, \rho)$, is entire,
and this conjecture is open in general. Again, one understands well
the {\it abelian} situation, i.e., when $\rho$ is a $1$-dimensional
representation; the kernel for such a $\rho$ defines an abelian
extension of $\Q$. By class field theory, such a $\rho$ is
associated to a character $\xi$ of finite order of the {idele class
group} $\A^\ast/\Q^\ast$; here, being {\it associated} means they
have the same $L$-function, with $L(s, \xi)$ being the one
introduced by Hecke, {\it albeit} in a different language. As usual,
we are denoting by $\A =\R\times\A_f$ the topological {\it ring of
adeles}, with $\A_f = \widehat{\Z}\otimes \Q$, and by $\A^\ast$ its
multiplicative group of {\it ideles}, which can be given the
structure of a locally compact abelian topological group with
discrete subgroup $\Q^\ast$.

Now fix a prime number $\ell$, and an algebraic closure $\overline
\Q_\ell$ of the field of $\ell$-adic numbers $\Q_\ell$, equipped
with an embedding $\overline \Q \hookrightarrow \overline
\Q_\ell$. Consider the set $\scR_\ell(n,\Q)$ of continuous,
semisimple representations
$$
\rho_\ell: \, \scG_\Q \, \rightarrow \, {\rm
GL}(n,\overline\Q_\ell),
$$
up to equivalence. The image of $\scG_\Q$ in such a representation
is usually not finite, and the simplest example of that is given
by the $\ell$-adic cyclotomic character $\chi_\ell$ given by the
action of $\scG_\Q$ on all the $\ell$-power roots of unity in
$\overline\Q$. Another example is given by the $2$-dimensional
$\ell$-adic representation on all the $\ell$-power {\it division
points} of an elliptic curve $E$ over $\Q$.

\medskip

The correct extension to the non-abelian case of the {\it idele
class character}, which appears in class field theory, is the notion
of an irreducible {\it automorphic representation} $\pi$ of GL$(n)$.
Such a $\pi$ is in particular a representation of the locally
compact group GL$(n,\A_F)$, which is a restricted direct product of
the local groups GL$(n, \Q_v)$, where $v$ runs over all the primes
$p$ and $\infty$ (with $\Q_\infty=\R$). There is a corresponding
factorization of $\pi$ as a tensor product $\otimes_v \pi_v$, with
all but a finite number of $\pi_p$ being {\it unramified}, i.e.,
admitting a vector fixed by the maximal compact subgroup $K_v$. At
the archimedean place $\infty$, $\pi_\infty$ corresponds to an
$n$-dimensional, semisimple representation $\sigma(\pi_\infty)$ of
the real Weil group $W_\R$, which is a non-trivial extension of
Gal$(\C/\R)$ by $\C^\ast$. Globally, by Schur's lemma, the center
$Z(\A)\simeq \A^\ast$ acts by a quasi-character $\omega$, which must
be trivial on $\Q^\ast$ by the automorphy of $\pi$, and so defines
an idele class character. Let us restrict to the central case when
$\pi$ is essentially unitary. Then there is a (unique) real number
$t$ such that the twisted representation $\pi_u:= \pi(t)=\pi\otimes
\vert\cdot\vert^t$ is unitary (with unitary central character
$\omega_u$). We are, by abuse of notation, writing
$\vert\cdot\vert^t$ to denote the quasi-character
$\vert\cdot\vert^t\circ{\rm det}$ of GL$(n,\A)$, where
$\vert\cdot\vert$ signifies the adelic absolute value, which is
trivial on $\Q^\ast$ by the Artin product formula. Roughly speaking,
to say that $\pi$ is automorphic means $\pi_u$ appears (in a weak
sense) in $L^2(Z(\A){\rm GL}(n,\Q)\backslash {\rm GL}(n,\A),
\omega_u)$, on which GL$(n,\A_F)$ acts by right translations. A
function $\varphi$ in this $L^2$-space whose averages over all the
horocycles are zero is called a {\it cusp form}, and $\pi$ is called
{\it cuspidal} if $\pi_u$ is generated by the right
GL$(n,\A_F)$-translates of such a $\varphi$. Among the automorphic
representations of GL$(n,\A)$ are certain distinguished ones called
{\it isobaric automorphic representations}. Any isobaric $\pi$ is of
the form $\pi_1\boxplus\pi_2\boxplus\dots\boxplus\pi_r$, where each
$\pi_j$ is a cuspidal representation of GL$(n_j,\A)$, such that
$(n_1,n_2,\dots,n_r)$ is a partition of $n$, where $\boxplus$
denotes the Langlands sum (coming from his theory of {\it Eisenstein
series}); moreover, every {\it constituent} $\pi_j$ is unique up to
isomorphism. Let $\scA(n,\Q)$ denote the set of isobaric automorphic
representations of GL$(n,\A)$ up to equivalence. Every isobaric
$\pi$ has an associated $L$-function $L(s,\pi)=\prod_v L(s,\pi_v)$,
which admits a meromorphic continuation and a functional equation.
Concretely, one associates at every prime $p$ where $\pi$ is
unramified, a conjugacy class $A(\pi)$ in GL$(n,\C)$, or
equivalently, an unordered $n$-tuple $(\alpha_{1,p}, \alpha_{2,p},
\dots, \alpha_{n,p})$ of complex numbers so that
$$
L(s,\pi_p) \, = \, \prod_{j=1}^n (1-\alpha_{j,p}p^{-s})^{-1}.
$$
If $\pi$ is cuspidal and non-trivial, $L(s,\pi)$ is entire; so is
the incomplete one $L^S(s,\pi)$ for any finite set $S$ of places of
$\Q$.

\medskip

Now suppose $\rho_\ell$ is an $n$-dimensional, semisimple
$\ell$-adic representation of $\scG_\Q={\rm Gal}(\overline\Q/\Q)$
{\it corresponds} to an automorphic representation $\pi$ of
GL$(n,\A)$. We will take this to mean that there is a finite set
$S$ of places including $\ell, \infty$ and all the primes where
$\rho_\ell$ or $\pi$ is ramified, such that we have,
$$
L(s,\pi_p) \, = \, L_p(s,\rho_\ell), \, \, \forall p\notin
S,\leqno(0.1)
$$
where the Galois Euler factor on the right is given by the
characteristic polynomial of $Fr_p$, the Frobenius at $p$, acting
on $\rho_\ell$. When (0.1) holds (for a suitable $S$), we will
write
$$
\rho_\ell \, \leftrightarrow \, \pi.
$$

A natural question in such a situation is to ask if $\pi$ is
cuspidal when $\rho_\ell$ is irreducible, and {\it vice-versa}. It
is certainly what is predicted by the general philosophy. However,
proving it is another matter altogether, and positive evidence is
scarce beyond $n=2$.

One can answer this question in the affirmative, for any $n$, if one
restricts to those $\rho_\ell$ which have {\it finite} image. In
this case, it also defines a continuous, $\C$-representation $\rho$,
the kind studied by E.~Artin ([A]). Indeed, the hypothesis implies
the identity of $L$-functions
$$
L^S(s,\rho \otimes \rho^\vee) \, = \, L^S(s, \pi \times
\pi^\vee),\leqno(0.2)
$$
where the superscript $S$ signifies the removal of the Euler factors
at places in $S$, and $\rho^\vee$ (resp. $\pi^\vee$) denotes the
contragredient of $\rho$ (resp. $\pi$). The $L$-function on the
right is the Rankin-Selberg $L$-function, whose mirific properties
have been established in the independent and complementary works of
Jacquet, Piatetski-Shapiro and Shalika ([JPSS], and of Shahidi
([Sh1,2]); see also [MW]. A theorem of Jacquet and Shalika ([JS1])
asserts that the {\it order of pole} at $s=1$ of $L^S(s,\pi \times
\pi^\vee)$ is $1$ iff $\pi$ is cuspidal. On the other hand, for any
finite-dimensional $\C$-representation $\tau$ of $\scG_\Q$, one has
$$
-{\rm ord}_{s=1}L^S(s,\tau) \, = \, {\rm dim}_\C {\rm
Hom}_{\scG_\Q}(\underline{1},\tau),\leqno(0.3)
$$
where $\underline{1}$ denotes the trivial representation of
$\scG_\Q$. Applying this with $\tau=\rho\otimes\rho^\vee \simeq {\rm
End}(\rho)$, we see that the order of pole of $L^S(s,
\rho\otimes\rho^\vee)$ at $s=1$ is $1$ iff the only operators in
End$(\rho)$ which commute with the $\scG_\Q$-action are scalars,
which means by Schur that $\rho$ is irreducible. Thus, {\it in the
Artin case}, {\it $\pi$ is cuspidal iff $\rho_\ell$ is irreducible}.

For general $\ell$-adic representations $\rho_\ell$ of $\scG_\Q$,
the order of pole at the right edge is not well understood. When
$\rho_\ell$ comes from {\it arithmetic geometry}, i.e., when it is a
Tate twist of a piece of the cohomology of a smooth projective
variety over $\Q$ which is cut out by algebraic projectors, an
important {\it conjecture of Tate} asserts an analogue of $(0.3)$
for $\tau=\rho_\ell\otimes\rho_\ell^\vee$, but this is unknown
except in a few families of examples, such as those coming from the
theory of {\it modular curves}, {\it Hilbert modular surfaces} and
{\it Picard modular surfaces}. So one has to find a different way to
approach the problem, which works at least in low dimensions.

\medskip

The main result of this paper is the following:

\medskip

\noindent{\bf Theorem A} \, \it Let $n \leq 5$ and let $\ell$ be a
prime. Suppose $\rho_\ell \, \leftrightarrow \, \pi$, for an
isobaric, algebraic automorphic representation $\pi$ of GL$(n, \A)$, and a
continuous, $\ell$-adic representation $\rho_\ell$ of $\scG_\Q$.
Assume
\begin{enumerate}
\item[(i)]$\rho_\ell$ is irreducible
\item[(ii)]$\pi$ is odd if $n\geq 3$
\item[(iii)]$\pi$ is semi-regular if $n=4$, and regular if $n=5$
\end{enumerate}
Then $\pi$ is cuspidal. \rm

\medskip

Some words of explanation are called for at this point. An isobaric
automorphic representation $\pi$ is said to be {\it algebraic}
([C$\ell$1]) if the restriction of $\sigma(\pi_\infty)$ to $\C^\ast$
is of the form $\oplus_{j=1}^n \chi_j$, with each $\chi_j$
algebraic, i.e., of the form $z\to z^{p_j}\overline{z}^{q_j}$ with
$p_j, q_j \in \Z$. (We do not assume that our automorphic
representations are unitary, and the arrow
$\pi_\infty\to\sigma(\pi_\infty)$ will be normalized
arithmetically.) For $n=1$, an algebraic $\pi$ is an idele class
character of type $A_0$ in the sense of Weil. One says that $\pi$ is
{\it regular} iff $\sigma(\pi_\infty)\vert_{\C^\ast}$ is a direct
sum of characters $\chi_j$, each occurring with {\it multiplicity
one}. And $\pi$ is {\it semi-regular} ([BHR]) if each $\chi_j$
occurs with {\it multiplicity at most two}. Suppose $\xi$ is a
$1$-dimensional representation of $W_\R$. Then, since $W_\R^{\rm ab}
\simeq \R^\ast$, $\xi$ is defined by a character of $\R^\ast$ of the
form $x\to \vert x\vert^w\cdot sgn(x)^{a(\xi)}$, with $a(\xi)\in
\{0,1\}$; here $sgn$ denotes the sign character of $\R^\ast$. For
every $w$, let $\sigma_\infty[\xi]:=
\sigma(\pi_\infty(\frac{1-n}{2}))[\xi]$ denote the isotypic
component of $\xi$, which has dimension at most $2$ (resp. $1$) if
$\pi$ is semi-regular (resp. regular), and is acted on by
$\R^\ast/\R_+^\ast \simeq \{\pm 1\}$. We will call a semi-regular
$\pi$ {\it odd} if for every character $\xi$ of $W_\R$, the
eigenvalues of $\R^\ast/\R_+^\ast$ on the $\xi$-isotypic component
are distinct. Clearly, any regular $\pi$ is odd under this
definition. See section $1$ for a definition of this concept for any
algebraic $\pi$, not necessarily semi-regular.

\medskip

I want to thank the organizers, Jae-Hyun-Yang in particular, and the
staff, of the {\it International Symposium on Representation Theory
and Automorphic Forms} in Seoul, Korea, first for inviting me to
speak there (during February $14-17, 2005$), and then for their
hospitality while I was there. The talk I gave at the conference was
on a different topic, however, and dealt with my ongoing work with
Dipendra Prasad on {\it selfdual representations}. I would also like
to thank F.~Shahidi for helpful conversations and the referee for
his comments on an earlier version, which led to an improvement of
the presentation. It is perhaps apt to end this introduction at this
point by acknowledging support from the National Science Foundation
via the grant $DMS-0402044$.

\bigskip

\section{Preliminaries}

\bigskip

\subsection{Galois Representations}

For any field $k$ with algebraic closure, denote by $\scG_k$ the
{\it absolute Galois group} of $\overline k$ over $k$. It is a
projective limit of the automorphism groups of finite Galois
extensions $E/k$. We furnish $\scG_k$ as usual with the {\it
profinite topology}, which makes it a {\it compact, totally
disconnected topological group}. When $k=\F_p$, there is for every
$n$ a unique extension of degree $n$, which is Galois, and
$\scG_{\F_p}$ is isomorphic to $\widehat{\Z} \simeq \lim_n \Z/n$,
topologically generated by the {\it Frobenius automorphism} $x\to
x^p$.

At each prime $p$, let $\scG_p$ denote the local Galois group
Gal$(\overline\Q_p/\Q_p)$ with inertia subgroup $I_p$, which fits
into the following exact sequence:
$$
1 \, \to \, I_p \, \to \scG_p \, \to \, \scG_{\F_p} \, \to \,
1.\leqno(1.1.1)
$$
The fixed field of $\overline \Q_p$ under $I_p$ is the {\it maximal
unramified extension} $\Q_p^{\rm ur}$ of $\Q_p$, which is generated
by all the roots of unity of order prime to $p$. One gets a natural
isomorphism of Gal$(\Q_p^{\rm ur}/\Q_p)$ with $\scG_{\F_p}$. If
$K/\Q$ is unramified at $p$, then one can lift the Frobenius element
to a conjugacy class $\varphi_p$ in Gal$(K/\Q)$.

All the Galois representations considered here will be continuous
and finite-dimensional. Typically, we will fix a prime $\ell$, and
algebraic closure $\overline \Q_\ell$ of the field $\Q_\ell$ of
$\ell$-adic numbers, and consider a continuous homomorphism
$$
\rho_\ell: \, \scG_\Q \, \rightarrow {\rm GL}(V_\ell),\leqno(1.1.2)
$$
where $V_\ell$ is an $n$-dimensional vector space over $\overline
\Q_\ell$. We will be interested only in those $\rho_\ell$ which are
unramified only at a finite set $S$ of primes. Then $\rho_\ell$
factors through a representation of the quotient group
$\scG_S:=G(\Q_S/\Q)$, where $\Q_S$ is the maximal extension of $\Q$
which is unramified out side $S$. One has the Frobenius classes
$\phi_p$ in $\scG_S$ for all $p\notin S$, and this allows one to
define the $L$-factors (with $s\in \C$)
$$
L_p(s, \rho_\ell) \, = \, {\rm det}\left(I-\varphi_pp^{-s}\, \vert
\, V_\ell\right)^{-1}.\leqno(1.1.3)
$$
Clearly, it is the reciprocal of a polynomial in $p^{-s}$ of degree
$n$, with constant term $1$, and it depends only on the equivalence
class of $\rho_\ell$. One sets
$$
L^S(s,\rho_\ell) \, = \, \prod_{p\notin S} \,
L_p(s,\rho_\ell).\leqno(1.1.4)
$$
When $\rho_\ell$ is the trivial representation, it is unramified
everywhere, and $L^S(s,\rho_\ell)$ is none other than the {\it Riemann
zeta function}. To define the {\it bad factors} at $p$ in
$S-\{\ell\}$, one replaces $V_\ell$ in this definition by the
subspace $V_\ell^{I_p}$ of {\it inertial invariants}, on which
$\varphi_p$ acts.

We are primarily interested in {\it semisimple representations} in
this article, which are direct sums of {\it simple} (or {\it
irreducible}) representations. Given any representation $\rho_\ell$
of $\scG_\Q$, there is an associated {\it semisimplification},
denoted $\rho_\ell^{\rm ss}$, which is a direct sum of the simple
Jordan-Holder components of $\rho_\ell$. A {\it theorem of
Tchebotarev} asserts the density of the Frobenius classes in the
Galois group, and since the local $p$-factors of $L(s,\rho_\ell)$
are defined in terms of the {\it inverse roots} of $\varphi_p$, one
gets the following standard, but useful result.

\medskip

\noindent{\bf Proposition 1.1.5} \, \it Let $\rho_\ell, \rho'_\ell$
be continuous, $n$-dimensional $\ell$-adic representations of
$\scG_\Q$. Then
$$
L^S(s,\rho_\ell) = L^S(s,\rho'_\ell) \, \implies \, \rho_\ell^{\rm
ss} \simeq {\rho'_\ell}^{\rm ss}.
$$
\rm

\medskip

The Galois representations $\rho_\ell$ which have {\it finite image}
are special, and one can view them as continuous
$\C$-representations $\rho$. Artin studied these deeply, and showed,
using the results of Brauer and Hecke, that the corresponding $L$-functions admit
meromorphic continuation and a functional equation of the form
$$
L^\ast(s,\rho) \, = \, \varepsilon(s,\rho)L^\ast(1-s,\rho^\vee),
\leqno(1.1.6)
$$
where $\rho^\vee$ denotes the contragredient representation on the
dual vector space, where
$$
L^\ast(s,\rho) \, = \, L(s,\rho)L_\infty(s,\rho),\leqno(1.1.7)
$$
with the {\it archimedean factor} $L_\infty(s,\rho)$ being a
suitable product (shifted) gamma functions. Moreover,
$$
\varepsilon(s,\rho) \, = \, W(\rho)N(\rho)^{s-1/2},\leqno(1.1.8)
$$
which is an entire function of $s$, with the (non-zero) $W(\rho)$
being called the {\it root number} of $\rho$. The scalar $N(\rho)$
is an integer, called the {\it Artin conductor} of $\rho$, and the
finite set $S$ which intervenes is the set of primes dividing
$N(\rho)$. The functional equation shows that
$W(\rho)W(\rho^\vee)=1$, and so $W(\rho)=\pm 1$ when $\rho$ is {\it
selfdual} (which means $\rho\simeq\rho^\vee$). Here is a useful
fact:

\medskip

\noindent{\bf Proposition 1.1.9} ([T]) \, \it Let $\tau$ be a
continuous, finite-dimensional $\C$-representation of $\scG_\Q$,
unramified outside $S$. Then we have
$$
-{\rm ord}_{s=1}L^S(s,\tau) \, = \, {\rm
Hom}_{\scG_\Q}(\underline{1}, \tau).
$$
\rm

\medskip

\noindent{\bf Corollary 1.1.10} \, Let $\rho$ be a continuous,
finite-dimensional $\C$-representation of $\scG_\Q$, unramified
outside $S$. Then $\rho$ is {\it irreducible} iff the incomplete
$L$-function $L^S(s,\rho\otimes\rho^\vee)$ has a {\it simple pole}
at $s=1$. \rm

\medskip

Indeed, if we set
$$
\tau: = \, \rho\otimes\rho^\vee \, \simeq \, {\rm
End}(\rho),\leqno(1.1.11)
$$
then Proposition $1.1.9$ says that the {\it order of pole} of
$L(s,\rho\otimes\rho^\vee)$ at $s=1$ is the {\it multiplicity of the
trivial representation} in End$(\rho)$ is $1$, i.e., iff the {\it
commutant} End$_{\scG_\Q}(\rho)$ is one-dimensional (over $\C$),
which in turn is equivalent, by Schur's lemma, to $\rho$ being
irreducible. Hence the Corollary.

\medskip

For general $\ell$-adic representations $\rho_\ell$, there is no
known analogue of Proposition $1.1.9$, though it is predicted to
hold (at the right edge of absolute convergence) by a {\it
conjecture of Tate} when $\rho_\ell$ comes from {\it arithmetic
geometry} (see [Ra4], section 1, for example). Tate's conjecture is
only known in certain special situations, such as for {\it $CM$
abelian varieties}. For the $L$-functions in Tate's set-up, say of
motivic weight $2m$, one does not even know that they make sense at
the {\it Tate point} $s=m+1$, let alone know its order of pole
there. Things get even harder if $\rho_\ell$ does not arise from a
geometric situation. One cannot work in too general a setting, and
at a minimum, one needs to require $\rho_\ell$ to have some good
properties, such as being unramified outside a finite set $S$ of
primes. Fontaine and Mazur conjecture ([FoM]) that $\rho_\ell$ is
{\it geometric} if it has this property (of being unramified outside
a finite $S$) and is in addition {\it potentially semistable}.

\medskip

\subsection{Automorphic Representations}

\medskip

Let $F$ be a number field with adele ring $\A_F=F_\infty\times\A_{F,f}$,
equipped with the adelic absolute value $\vert\cdot\vert=\vert\cdot\vert_\A$.
For every algebraic group $G$ over $F,$ let
$G(\A_F) = G(F_\infty)\times G(\A_{F,f})$ denote the
restricted direct product $\prod'_v G(F_v),$ endowed with the usual
locally compact topology.
Then $G(F)$ embeds in $G(\A_F)$ as a disct=rete subgroup,
and if $Z_n$ denotes the center of GL$(n)$, the homogeneous space
${\rm GL}(n,F)Z(\A_F)\backslash {\rm GL}(n,\A_F)$ has finite volume relative to the
relatively invariant quotient measure
induced by a Haar measure on ${\rm GL}(n,\A_F)$. An irreducible representation
$\pi$ of ${\rm GL}(n,\A_F)$ is admissible if it admits a
factorization as a restricted tensor product $\otimes'_v \pi_v$,
where each $\pi_v$ is admissible and
for almost all finite places $v$, $\pi_v$ is {\it unramified}, i.e., has a no-zero vector fixed by
$K_v= {\rm GL}(n,\scO_v)$. (Here, as usual, $\scO_v$ denotes the ring of integers of the local completion
$F_v$ of $F$ at $v$.)

Fixing a unitary idele class character $\omega$,
which can be viewed as a character of $Z_n(\A_F)$, we may consider the space
$$
L^2(n,\omega): = L^2({\rm GL}(n,F)Z(\A_F)\backslash {\rm GL}(n,\A_F), \omega),\leqno(1.2.1)
$$
which consists of (classes of) functions on GL$(n,\A_F)$ which are left-invariant under GL$(n,F)$,
transform under $Z(\A_F)$ according to $\omega$, and are square-integrable modulo
GL$(n,F)Z(\A_F)$. Clearly, $L^2(n,\omega)$ is a unitary representation of GL$(n,\A_F)$ under
the right translation action on functions.
The {\bf space of cusp forms}, denoted $L^2_0(n, \omega)$, consists of
functions $\varphi$ in $L^2(n,\omega)$ which satisfy the following for {\it every} unipotent
radical $U$ of a standard parabolic subgroup $P=MU$:
$$
\int_{U(F)\backslash U(\A_F)} \varphi(ux) \, = \, 0. \leqno(1.2.2)
$$
To say that $P$ is a standard parabolic means that it contains the {\it Borel subgroup}
of upper triangular matrices in GL$(n)$. A basic fact asserts that $L^2_0(n,\omega)$ is a
subspace of the discrete spectrum of $L^2(n,\omega)$.

By a {\bf unitary cuspidal} (automorphic) representation $\pi$ of ${\rm GL}_n(\A_F)$,
we will mean an irreducible, unitary representation occurring in
$L^2_0(n,\omega)$. We will, by abuse of notation, also denote the underlying admissible
representation by $\pi$. (To be precise, the unitary representation is on the Hilbert space
completion of the admissible space.) Roughly speaking, unitary automorphic representations
of GL$(n,\A_F)$ are
those which appear weakly in $L^2(n,\omega)$ for some $\omega$. We will refrain from
recalling the definition precisely, because we will work totally
with the subclass of
{\it isobaric automorphic representations}, for which one can
take Theorem 1.2.10 (of Langlands)
below as their definition.

\medskip

If $\pi$ is an admissible representation of GL$(n,\A_F)$, then for
any $z\in\C$, we define the {\it analytic Tate twist} of $\pi$ by
$z$ to be
$$
\pi(z): = \, \pi\otimes\vert\cdot\vert^{z},\leqno(1.2.3)
$$
where $\vert\cdot\vert^z$ denotes the $1$-dimensional representation
of GL$(n,\A_F)$ given by
$$
g \to \, \vert{\rm det}(g)\vert^z \, = \, e^{z\log(\vert{\rm
det}(g)\vert)}.
$$
Since the adelic absolute value $\vert\cdot\vert$ takes det$(g)$ to a
positive real number,
its logarithm is well defined.

In general, by a {\it cuspidal automorphic representation}, we will
mean an irreducible admissible representation of ${\rm GL}(n,\A_F)$
for which there exists a real number $w$, which we will call the
{\it weight of $\pi$} such that the Tate twist
$$
\pi_u:= \pi(w/2)\leqno(1.2.4)
$$
is a unitary cuspidal representation. Note that the central character of $\pi$
and of its unitary {\it avatar} $\pi_u$ are related as follows:
$$
\omega_{\pi} \, = \, \omega_{\pi_u}\vert\cdot\vert^{-nw/2},\leqno(1.2.5)
$$
which is easily checked by looking at the situation at the unramified
primes, which suffices.

\medskip

For any irreducible, automorphic representation $\pi$ of
$GL(n,\A_F),$ there is an associated $L$-function
$L(s, \pi) = L(s, \pi_{\infty})L(s, \pi_f)$, called the
{\it standard} $L-$function ([J]) of $\pi;$. It has an
Euler product expansion
$$
L(s,\pi) \, = \, \prod_v \, L(s, \pi_v), \leqno (1.2.6)
$$
convergent in a right-half plane. If $v$ is an archimedean place,
then one knows (cf. [La1]) how to associate a semisimple
$n-$dimensional $\C-$representation $\sigma(\pi_v)$ of the Weil
group $W_{F_v},$ and $L(\pi_v,s)$ identifies with $L(\sigma_v,s).$
We will {\it normalize} this correspondence $\pi_v \to \sigma(\pi_v)$ in such a way
that it respects algebraicity.
Moreover, if $v$ is a finite place where $\pi_v$ is
unramified, there is a corresponding semisimple
conjugacy class $A_v(\pi)$ (or $A(\pi_v)$) in GL$(n,\C)$ such that
$$
L(s,\pi_v) \, = \, {\rm {det}}(1-A_v(\pi)T)^{-1}|_{T=q_v^{-s}}.
\leqno (1.2.7)
$$
We may find a diagonal representative diag$(\alpha_{1,v}(\pi), ... ,
\alpha_{n,v}(\pi))$ for $A_v(\pi),$ which is unique up to
permutation of the diagonal entries. Let $[\alpha_{1,v}(\pi), ... ,
\alpha_{n,v}(\pi) ]$ denote the unordered $n-$tuple of complex numbers
representing $A_v(\pi)$. Since
$W_{F,v}^{{\rm {ab}}} \simeq F_v^\ast,$ $A_v(\pi)$ clearly defines
an abelian $n-$dimensional representation $\sigma(\pi_v)$ of
$W_{F,v}.$ If $\underline{1}$ denotes the trivial representation of
GL$(1,\A_F)$, which is cuspidal, we have
$$
L(s,\underline{1}) \, = \, \zeta_F(s),
$$
the Dedekind zeta function of $F$. (Strictly speaking, we should take
$L(s,\underline{1}_f)$ on the left, since the right hand side is missing
the archimedean factor, but this is not serious.)

The fundamental work of Godement and Jacquet, when used in conjunction
with the Rankin-Selberg theory (see $1.3$ below), yields the following:

\medskip

\noindent{\bf Theorem 1.2.8} ([J]) \quad \it Let $n \geq 1,$
and $\pi$ a non-trivial cuspidal automorphic representation of
GL$(n,\A_F)$. Then $L(s,\pi)$ is entire. Moreover, for any finite set $S$ of
places of $F,$ the incomplete $L-$function
$$L^S(s,\pi) \, = \,
\prod_{v \notin S} L(s,\pi_v)
$$
is holomorphic and non-zero in $\Re(s)> w+1$ if $\pi$ has weight
$w$. Moreover, there is a functional equation
$$
L(w+1-s,\pi^\vee) \, = \, \varepsilon(s,\pi)L(s,\pi)\leqno(1.2.9)
$$
with
$$
\varepsilon(s,\pi) \, = \, W(\pi)N_\pi^{(w+1)/2-s}.
$$
Here $N_\pi$ denotes the norm of the conductor ${\mathcal N}_\pi$ of $\pi$,
and $W(\pi)$ is the root number of $\pi$.
\rm

\medskip

Of course when $w=0$, i.e., when $\pi$ is unitary,
the statement comes to a more familiar form.
When $n = 1,$ a $\pi$ is simply an idele class character
and this result is due to Hecke.

\medskip

By the theory of Eisenstein series, there is a sum operation
$\boxplus$ ([La2], [JS1]):

\medskip

\noindent{\bf Theorem 1.2.10} ([JS1]) \quad \it Given any $m-$tuple
of cuspidal automorphic representations $\pi_1, ..., \pi_m$ of
GL$(n_1,\A_F), ... ,$ GL$(n_m,\A_F)$ respectively, there exists an
irreducible, automorphic representation $\pi_1 \boxplus ... \boxplus
\pi_m$ of GL$(n,\A_F),$ $n \, = \, n_1 + ... + n_m,$ which is unique
up to equivalence, such that for any finite set $S$ of places,
$$
L^S(s, \boxplus_{j=1}^m \pi_j) \, = \, \prod_{j=1}^m L^S(s, \pi_j).
\leqno(1.2.11)
$$
\rm

\medskip

Call such a (Langlands) sum $\pi \simeq \boxplus_{j=1}^m \pi_j$,
with each $\pi_j$ cuspidal, an {\it isobaric automorphic}, or just
{\it isobaric} (if the context is clear), representation. Denote by
ram$(\pi)$ the finite set of finite places where $\pi$ is ramified,
and let $\mathfrak N(\pi)$ be its conductor.

For every integer $n \geq 1,$ set:
$$
\mathcal  A(n,F) \, = \, \{\pi: {\rm {isobaric}} \, {\rm {
representation}} \, {\rm {of}} \, {\rm {GL}}(n,\A_F) \}/{\simeq },
\leqno (1.2.12)
$$
\noindent and
$$
\mathcal  A_0(n,F) \, = \, \{\pi \in \mathcal  A(n,F) | \, \pi \,
{\rm {cuspidal}} \}.
$$
Put $\mathcal  A(F) \, = \, \cup_{n \geq 1} \mathcal  A(n,F)$ and
$\mathcal  A_0(F) \, = \, \cup_{n \geq 1} \mathcal  A_0(n,F).$

\medskip

\noindent{\bf Remark 1.2.13}. \quad One can also define the analogs
of $\mathcal A(n, F)$ for local fields $F$, where the ``cuspidal''
subset $\mathcal A_0(n, F)$ consists of essentially
square-integrable representations of GL$(n, F)$. See [La3] (or
[Ra1]) for details.

\medskip

Given any polynomial representation
$$
r: \, {\rm GL}(n,\C) \, \rightarrow \, {\rm GL}(N,\C),\leqno(1.2.14)
$$
one can associate an $L$-function to the pair $(\pi,r)$, for any
isobaric automorphic representation $\pi$ of GL$(n,\A_F)$:
$$
L(s,\pi; r) \, = \, \prod_v \, L(s,\pi_v;r), \leqno(1.2.15)
$$
in such a way that at any finite place $v$ where $\pi$ is {\it unramified}
with residue field $\F_q$,
$$
L(s,\pi_v; r) \, = \, {\rm {det}}(1-A_v(\pi; r)T)^{-1}|_{T=q_v^{-s}},\leqno(1.2.16)
$$
with
$$
A_v(\pi; r) \, = \, r(A_v(\pi)).\leqno(1.2.17)
$$
The conjugacy class $A_v(\pi; r)$ in GL$(N,\C)$ is again represented by an
unordered $N$-tuple of complex numbers.

\medskip

The {\bf Principle of Functoriality} predicts the existence of an isobaric automorphic
representation $r(\pi)$ of GL$(N,\A_F)$ such that
$$
L(s, r(\pi)) \, = \, L(s,\pi; r)\leqno(1.2.18)
$$
A weaker form of the conjecture, which suffices for questions like what we are considering,
asserts that this identity holds outside a finite set $S$ of places.

This conjecture is known in the following cases of ${\bf (n,r)}$:

\noindent{${\bf (1.2.19)}$}
\begin{enumerate}
\item[]\quad ${\bf (2,{\rm sym}^2)}$: \, Gelbart-Jacquet ([GJ])
\item[]\quad ${\bf (2,{\rm sym}^3)}$: \, Kim-Shahidi ([KSh1])
\item[]\quad ${\bf (2,{\rm sym}^4)}$: \, Kim ([K])
\item[]\quad ${\bf (4,\Lambda^2)}$: \, Kim ([K])
\end{enumerate}

In this paper we will make use of the last instance of functoriality,
namely the {\it exterior square}
transfer of automorphic forms from GL$(4)$ to GL$(6)$.

\bigskip

\subsection{Rankin-Selberg $L$-functions}

\medskip

The results here are due to the independent and partly
complementary, deep works of Jacquet, Piatetski-Shapiro and Shalika,
and of Shahidi. Let $\pi,$ $\pi'$ be isobaric automorphic
representations in $\mathcal A(n,F),$ $\mathcal A(n',F)$
respectively. Then there exists an associated Euler product $L(s,
\pi \times \pi')$ ([JPSS], [JS1], [Sh1,2], [MW], [CoPS2]), which
converges in $\{\Re(s) > 1 \},$ and admits a meromorphic
continuation to the whole $s-$plane and satisfies the functional
equation, which is given in the unitary case by
$$
L(s, \pi \times \pi') \, = \, \varepsilon(s, \pi \times \pi')L(1-s,
\pi^\vee \times {\pi'}^\vee), \leqno(1.3.1)
$$
with
$$
\varepsilon(s, \pi \times \pi') \, = \, W(\pi \times \pi') N(\pi
\times \pi')^{\frac{1}{2}-s},
$$
where the {\sl conductor} $N(\pi \times \pi')$ is a positive integer
not divisible by any rational prime not intersecting the
ramification loci of $F/\Q$, $\pi$ and $\pi'$, while $W(\pi \times
\pi')$ is the {\it root number} in $\C^\ast$. As in the Galois case,
$W(\pi \times \pi')W(\pi^\vee \times {\pi'}^\vee) = 1$, so that
$W(\pi \times \pi') = \pm 1$ when $\pi, \pi'$ are self-dual.

It is easy to deduce the functional equation when $\pi, \pi'$ are
not unitary. If they are cuspidal of weights $w, w'$ respectively,
the functional equation relates $s$ to $w+w'+1-s$. Moreover, since
$\pi^\vee, {\pi'}^\vee$ have respective weights $-w,-w'$,
$\pi\times\pi^vee$ and $\pi'\times{\pi'}^\vee$ still have weight
$0$.

When $v$ is archimedean or a finite place unramified for $\pi,
\pi'$,
$$
L_v(s, \pi \times \pi') \, = \, L(s, \sigma(\pi_v) \otimes
\sigma(\pi'_v)). \leqno (1.3.2)
$$
In the archimedean situation, $\pi_v \to \sigma(\pi_v)$ is the arrow
to the representations of the Weil group $W_{F_v}$ given by [La1].
When $v$ is an unramified finite place, $\sigma(\pi_v)$ is defined
in the obvious way as the sum of one dimensional representations
defined by the Langlands class $A(\pi_v)$.

When $n=1,$ $L(s, \pi \times \pi') \, = \, L(s, \pi\pi'),$ and when
$n=2$ and $F = \Q,$ this function is the usual Rankin-Selberg
$L-$function, extended to arbitrary global fields by Jacquet.

\medskip

\noindent{\bf Theorem 1.3.3} ([JS1], [JPSS]) \quad \it Let $\pi \in
\mathfrak A_0(n,F)$, $\pi' \in \mathfrak A_0(n',F)$, and $S$ a
finite set of places. Then $L^S(s, \pi \times \pi')$ is entire
unless $\pi$ is of the form ${\pi'}^\vee \otimes |.|^w$, in which
case it is holomorphic outside $s = w, 1-w$, where it has simple
poles. \rm

\medskip

The {\it Principle of Functoriality} implies in this situation that
given $\pi, \pi'$ as above, there exists an isobaric automorphic
representation $\pi\boxtimes\pi'$ of GL$(nn',\A_F)$ such that
$$
L(s,\pi\boxtimes\pi') \, =, \, L(s,\pi\times\pi').\leqno(1.3.4)
$$

\medskip

The (conjectural) functorial product$\boxtimes$ is the automorphic
analogue of the usual tensor product of Galois representations. For
the importance of this product, see [Ra1], for example.

One can always construct $\pi\boxtimes\pi'$ as an {\it admissible}
representation of GL$(nn',\A_f)$, but the subtlety lies in showing
that this product is automorphic.

The automorphy of $\boxtimes$ is known in the following cases, which
will be useful to us:

\noindent{${\bf (1.3.5)}$}
\begin{enumerate}
\item[]\quad ${\bf (n,n')= (2,2)}$: \, ([Ra2])
\item[]\quad ${\bf (n,n')=(2,3)}$: \, Kim-Shahidi ([KSh1])
\end{enumerate}

\medskip

The reader is referred to section $11$ of [Ra4], which contains some
refinements, explanations, refinements and (minor) errata for [Ra2].
It may be worthwhile remarking that Kim and Shahidi use the
functorial product on GL$(2)\times $GL$(3)$ which they construct to
prove the {\it symmetric cube} lifting for GL$(2)$ mentioned in the
previous section (see $(1.2.11)$. A {\it cuspidality criterion} for
the image under this transfer is proved in [Ra-W], with an
application to the cuspidal cohomology of congruence subgroups of
SL$(6,\Z)$.

\bigskip

\subsection{Modularity and the problem at hand}

\medskip

The general Langlands philosophy asserts that if $\rho_\ell$ is an
$n$-dimensional $\ell$-adic representation of $\scG_\Q$ arising as a
factor of the cohomology of a smooth, projective variety over $\Q$,
then there is an isobaric automorphic representation $\pi$ of GL$(n,
\A)$ such that for a suitable finite set $S$ of places (including
$\infty$), we have an identity of the form
$$
L^S(s,\rho_\ell) \, = \, L^S(s,\pi).\leqno(1.4.1)
$$
When this happens, we will say that $\rho_\ell$ is {\it modular},
and we will write
$$
\rho_\ell \, \leftrightarrow \, \pi.\leqno(1.4.2)
$$
One says that $\rho_\ell$ is {\it strongly modular} if the identity
$(1.4.1)$ holds for the full $L$-function, i.e., with $S$ empty.

Recall from (the end of) section 1.1 that a striking conjecture of
Fontaine and Mazur ([FoM]) asserts that a Galois representation
$\rho_\ell$ comes from arithmetical geometry, as required in the
modularity conjecture above, if it is potentially semistable at
$\ell$ and has good reduction at almost all primes.

Special cases of the modularity conjecture were known earlier, the
most famous one being the modularity conjecture for the $\ell$-adic
representations $\rho_\ell$ defined by the Galois action on the
$\ell$-power division points of elliptic curves $E$ over $\Q$,
proved recently in the spectacular works of Wiles, Taylor, Diamond,
Conrad and Breuil.

\medskip

We will not consider any such (extremely) difficult question in this
article. Instead we will be interested in the following:

\medskip

\noindent{\bf Question 1.4.3} \, \it When a modular $\rho_\ell$ is
irreducible, is the corresponding $\pi$ cuspidal? And conversely?\rm

\medskip

This seemingly reasonable question turns out to be hard to check in
dimensions $n > 2$.

\medskip

One thing that is clear is that the $\pi$ associated to any
$\rho_\ell$ needs to be {\it algebraic} in the sense of Clozel
([C$\ell$1]). To define the notion of {\it algebraicity}, first
recall that by Langlands, the archimedean component $\pi_\infty$ is
associated to an $n$-dimensional representation
$\sigma(\pi_\infty)$, sometimes written $\sigma_\infty(\pi)$, of the
real Weil group $W_\R$, with corresponding equality of the
archimedean $L$-factors $L_\infty(s,\rho_\ell)$ and
$L(s,\pi_\infty)$. We will normalize things so that the
correspondence is algebraic. One can explicitly describe $W_\R$ as
$\C^\ast\cup j\C^\ast$, with $jzj^{-1}=\overline{z}$ and $j^2=-1$.
One gets a canonical exact sequence
$$
1 \, \to \, \C^\ast \, \to \, W_\R \, \to \, \scG_\R \, \to \,
1\leqno(1.4.4)
$$
which represents the unique non-trivial extension of $\scG_\R$ by
$\C^\ast$. One has a decomposition
$$
\sigma(\pi_\infty)\vert_{\C^\ast} \, \simeq \, \oplus_{j=1}^n \,
\xi_j,\leqno(1.4.5)
$$
where each $\xi_j$ is a (quasi-)character of $\C^\ast$. One says
that $\pi$ is {\it algebraic} when every one of the characters
$\chi_j$ is algebraic, i.e., there are integers $p_j, q_j$ such that
$$
\chi_j(z) \, = \, z^{p_j}\overline{z}^{q_j}.\leqno(1.4.6)
$$
This analogous to having a {\it Hodge structure}, which is what one
would expect if $\pi$ were to be related to a geometric object.

One says that $\pi$ {\it regular} if for all $i\ne j$,
$\chi_i\ne\chi_j$. In other words, each character $\chi_j$ appears
in the restriction of $\sigma_\infty(\pi)$ to $\C^\ast$ with {\it
multiplicity one}. We say (following [BHR]) that $\pi$ is {\it
semi-regular} if the multiplicity of each $\chi_j$ is at most $2$.

When $n=2$, any $\pi$ defined by a classical {\it holomorphic
newform} $f$ of weight $k\geq 1$ is algebraic and semi-regular. It
is regular iff $k\geq 2$. One also expects any {\it Maass waveform}
$\varphi$ of {\it weight $0$} and {\it eigenvalue $1/4$} for the
{\it hyperbolic Laplacian} to be algebraic; there are interesting
examples of this sort coming from the work of Langlands (resp.
Tunnell) on tetrahedral (resp. octahedral) Galois representations
$\rho$ which are {\it even}; the {\it odd} ones correspond to
holomorphic newforms of weight $1$. We will not consider the even
situation in this article.

\medskip

Given a holomorphic newform $f(z)= \sum_{n=1}^\infty a_nq^n$,
$q=e^{2\pi iz}$, of weight $2$, resp. $k \geq 3$, resp. $k=1$, level
$N$ and character $\omega$, one knows by Eichler and Shimura, resp.
Deligne ([De]), resp Deligne-Serre ([DeS]), that there is a
continuous, irreducible representation
$$
\rho_\ell: \scG_\Q \, \rightarrow \, {\rm
GL}(2,\Q_\ell)\leqno(1.4.7)
$$
such that for all primes $p\nmid N\ell$,
$$
{\rm tr}\left( Fr_p \, \vert \, \rho_\ell\right) \, = \, a_p
$$
and
$$
{\rm det}(\rho_\ell) \, = \, \omega\chi_\ell^{k-1},
$$
where $\chi_\ell$ is the $\ell$-adic {\it cyclotomic character}  of
$\scG_\Q$, given by the Galois action on the $\ell$-power roots of
unity in $\overline \Q$, and $Fr_p$ is the {\it geometric Frobenius}
at $p$, which is the inverse of the arithmetic Frobenius.

\bigskip

\subsection{Parity}

\medskip

We will first first introduce this crucial concept over the base
field $\Q$, as that is what is needed in the remainder of the
article.

We will need to restrict our attention to those isobaric forms $\pi$
on GL$(n)/\Q$ which are odd in a suitable sense. It is instructive
to first consider the case of a classical holomorphic newform $f$ of
weight $k\geq 1$ and character $\omega$ relative to the congruence
subgroup $\Gamma_0(N)$. Since $\Gamma_0(N)$ contains $-I$, it
follows that $\omega(-1) = (-1)^k$. One could be tempted to call $a
\pi$ defined by such an $f$ to be even (or odd) according as
$\omega$ is even (or odd), but it would be a wrong move. One should
look not just at $\omega$, but at the determinant of the associated
$\rho_\ell$, i.e., the $\ell$-adic character
$\omega\chi_\ell^{k-1}$, which is odd for all $k$ ! So all such
$\pi$ defined by holomorphic newforms are {\it arithmetically} odd.
The only even ones for GL$(2)$ are (analytic Tate twists of) Maass
forms of weight $0$ and Laplacian eigenvalue $1/4$.

\medskip

The maximal abelian quotient of $W_\R$ is $\R^\ast$, and the
restriction of the abelianization map to $\C^\ast$ identifies with
the norm map $z \to \vert z\vert$. So every (quasi)-character $\xi$
of $W_\R$ identifies with one of $\R^\ast$, given by $x\to
sgn(x)^a\vert x\vert^t$ for some $t$, with $a \in \{0,1\}$. Clearly,
$\xi$ determines, and is determined by $(t,a)$. If $\pi$ is an
isobaric automorphic representation, let $\sigma_\infty[\xi]$
denote, for each such $\xi$, the $\xi$-isotypic component of
$\sigma_\infty(\pi)$. The {\it sign group} $\R^\ast/\R_+^\ast$ acts
on each isotypic component. Let $m_+(\pi,\xi)$ (resp.
$m_-(\pi,\xi)$) denote the multiplicity of the eigenvalue $+1$
(resp. $-1$), under the action of $\R^\ast/\R_+^\ast$ on
$\sigma_\infty(\xi)$.

\medskip

\noindent{\bf Definition 1.5.1} \, \it Call an isobaric automorphic
representation $\pi$ of GL$(n,\A)$ \emph{odd} if for every
one-dimensional representation $\xi$ of $W_\R$ occurring in
$\sigma_\infty(\pi)$,
$$
\vert m_+(\pi,\xi) - m_-(\pi,\xi)\vert \, \leq \, 1.
$$
\rm

\medskip

Clearly, when the dimension of $\sigma_\infty[\xi]$ is even, the
multiplicity of $+1$ as an eigenvalue of the sign group needs to be
equal to the multiplicity of $-1$ as an eigenvalue.

Under this definition, all forms on GL$(1)/\Q$ are odd. So are the
$\pi$ on GL$(2)/\Q$ which are defined by holomorphic newforms of
weight $k \geq 2$. The reason is that $\pi_\infty$ is (for $k\geq
2$) a discrete series representation, and the corresponding
$\sigma_\infty(\pi)$ is an irreducible $2$-dimensional
representation of $W_\R$ induced by the (quasi)-character $z\to
z^{-(k-1)}$ of the subgroup $\C^\ast$ of index $2$, and our
condition is vacuous. On the other hand, if $k=1$,
$\sigma_\infty(\pi)$ is a reducible $2$-dimensional representation,
given by $\underline{1} \oplus sgn$. The eigenvalues are $1$ on
$\sigma_\infty(\underline{1})$ and $-1$ on $\sigma_\infty(sgn)$. On
the other hand, a Maass form of weight $0$ and $\lambda=1/4$, the
eigenvalue $1$ (or $-1$) occurs with multiplicity $2$, making the
$\pi$ it defines an even representation. So our definition is a good
one and gives what we know for $n=2$.

\medskip

For any $n$, note that if $\pi$ is algebraic and regular, it is
automatically odd. If $\pi$ is algebraic and semi-regular, each
isotypic space is one or two-dimensional, and in the latter case, we
want both eigenvalues to occur for $\pi$ to be odd.

\medskip

Finally, if $F$ is any number field with a real place $u$, we can
define, in exactly the same way, when an algebraic, isobaric
automorphic representation of GL$(n,\A_F)$ is {\it arithmetically
odd at} ${\it u}$. If $F$ is {\it totally real}, then we say that
$\pi$ is {\it totally odd} if it is so at {\it every} archimedean
place.

\bigskip

\section{The first step in the proof}

\bigskip

Let $\rho_\ell, \pi$ be as in Theorem A. Since $\rho_\ell$ is
irreducible, it is in particular semisimple. Suppose $\pi$ is not
cuspidal. We will obtain a contradiction.

\medskip

\noindent{\bf Proposition 2.1} \, \it Let $\rho_\ell, \pi$ be
associated, with $\pi$ algebraic, semi-regular and odd. Suppose we
have, for some $r
> 1$, an isobaric sum decomposition
$$
\pi \, \simeq \, \boxplus_{j=1}^r \, \eta_j,\leqno(2.2)
$$
where each $\eta_j$ a cuspidal automorphic representation of
GL$(n_j,\A)$, with $n_j\leq 2$  ($\forall j)$. Then $\rho_\ell$
cannot be irreducible.\rm

\medskip

\noindent{\bf Corollary 2.3} \, \it Theorem A holds when $\pi$
admits an isobaric sum decomposition such as $(2.2)$ with each
$n_j\leq 2$. In particular, it holds for $n\leq 3$. \rm

\medskip

\noindent{\bf Proof of Proposition}. \, The hypothesis that $\pi$ is
algebraic and semi-regular implies easily that each $\eta_j$ is also
algebraic and semi-regular. Let $J_m$ denote the set of $j$ where
$n_j=m$.

First look at any $j$ in $J_1$. Then the corresponding $\eta_j$ is
an idele class character. Its algebraicity implies that, in
classical terms, it corresponds to an algebraic Hecke character
$\nu_j$. By Serre ([Se]), we may attach an abelian $\ell$-adic
representation $\nu_{j,\ell}$ of $\scG_\Q$ of dimension $1$. It
follows that for some finite set $S$ of places containing $\ell$,,
$$
L^S(s,\nu_{j,\ell}) \, = \, L^S(s,\eta_j) \quad {\rm whenever}
\quad n_j=1.\leqno(2.4)
$$

Next consider any $j$ in $J_2$. If $\sigma(\eta_{j,\infty})$ is
irreducible, then a twist of $\eta_j$ must correspond to a classical
holomorphic newform $f$ of weight $k\geq 2$. Moreover, the
algebraicity of $\eta_j$ forces this twist to be algebraic. Hence by
Deligne, there is a continuous representation
$$
\tau_{j,\ell}: \scG_\Q \, \rightarrow \, {\rm
GL}(2,\overline\Q_\ell),\leqno(2.5)
$$
ramified only at a finite of primes such that at every $p\ne \ell$
where the representation is unramified,
$$
{\rm tr}(Fr_p \, \vert \, \tau_{j,\ell}) \, = \,
a_p(\eta_j),\leqno(2.6)
$$
and the determinant of $\tau_{j,\ell}$ corresponds to the central
character $\omega_j$ of $\eta_j$. Moreover, $\tau_{j,\ell}$ is
irreducible, which is not crucial to us here.

We also need to consider the situation, for any fixed $j \in J_2$,
when $\sigma(\eta_{j,\infty})$ is reducible, say of the form
$\chi_1\oplus\chi_2$. Since $\eta_j$ is cuspidal, by the {\it
archimedean purity} result of Clozel ([C$\ell 1$]),
$\chi_1\chi_2^{-1}$ must be $1$ or $sgn$. The former cannot happen
due to the oddness of $\pi$. It follows that $\eta_j$ is defined by
a classical holomorphic newform $f$ of weight $1$, and by a result
of Deligne and Serre ([DeS]), there is a $2$-dimensional $\ell$-adic
representation $\tau_{j,\ell}$ of $\scG_\Q$ {\it with finite image},
which is irreducible, such that $(2.6)$ holds.

Since the set of Frobenius classes $Fr_p$, as $p$ runs over primes
outside $S$, is dense in the Galois group by Tchebotarev, we must
have, by putting all these cases together,
$$
\rho_\ell \, \simeq \, \left(\oplus_{j\in J_1} \,
\nu_{j,\ell}\right) \, \oplus \, \left(\oplus_{j\in J_2} \,
\tau_{j,\ell}\right),\leqno(2.7)
$$
which contradicts the irreducibility of $\rho_\ell$, since by
hypothesis, $r=\vert J_1\vert +\vert J_2\vert \, \geq 2$.

\qed

\bigskip

\section{The second step in the proof}

\bigskip

Let $\rho_\ell, \pi$ be as in Theorem A. Suppose $\pi$ is not
cuspidal. In view of Proposition 2.1, we need only consider the
situation where $\pi$ is an isobaric sum $\boxplus_j \eta_j$, with
an $\eta_j$ being a cusp form on GL$(m)/\Q$ for some $m\geq 3$.

\medskip

\noindent{\bf Proposition 3.1} \, \it Let $\rho_\ell, \pi$ be
associated, with $\pi$ an algebraic cusp form on GL$(n)/\Q$ which is
semi-regular and odd. Suppose we have an isobaric sum decomposition
$$
\pi \, \simeq \, \eta \boxplus \eta',\leqno(3.2)
$$
where $\eta$ is a cusp form on GL$(3)/\Q$ and $\eta'$ an isobaric
automorphic representation of GL$(r, \A)$ for some $r\geq 1$.
Moreover, assume that there is an $r$-dimensional $\ell$-adic
representation $\tau_\ell'$ of $\scG_\Q$ associated to $\eta'$.
Then we have the isomorphism of $\scG_\Q$-modules:
$$
\rho_\ell^\vee \oplus (\rho_\ell\otimes\tau'_\ell) \, \simeq \,
\Lambda^2(\rho_\ell)\oplus{\tau'_\ell}^\vee\oplus{\rm
sym}^2(\tau'_\ell).\leqno(3.3)
$$
\rm

\medskip

\noindent{\bf Corollary 3.4} \, \it Let $\rho_\ell, \pi$ be
associated, with $\pi$ algebraic, semi-regular and odd. Suppose
$\pi$ admits an isobaric sum decomposition such as $(3.2)$ with
$r\leq 2$. Then $\rho_\ell$ is reducible. \rm

\medskip

\noindent{\bf Proposition $3.1 \implies$ Corollary $3.4$}: When
$r\leq 2$, $\eta'$ is either an isobaric sum of algebraic Hecke
characters or cuspidal, in which case, thanks to the oddness, it is
defined by a classical cusp form on GL$(r)/\Q$ of weight $\geq 1$.
In either case we have, as seen in the previous section, the
existence of the associated $\ell$-adic representation $\tau'_\ell$,
which is irreducible exactly when $\eta_\ell$ is cuspidal. Then by
the Proposition, the decomposition $(3.3)$ holds. If $r=1$ or $r=2$
with $\eta'$ Eisensteinian, $(3.3)$ implies that a $1$-dimensional
representation (occurring in $\tau'_\ell$) is a summand of a twist
of either $\rho_\ell$ or $\rho_\ell^\vee$. Hence the Corollary.

\qed

\medskip

Combining Corollary $2.3$ and Corollary $3.4$, we see that the
irreducibility of $\rho_\ell$ forces the corresponding $\pi$ to be
cuspidal when $n\leq 4$ under the hypotheses of Theorem A. So we
obtain the following:

\medskip

\noindent{\bf Corollary 3.5} \, \it Theorem A holds for $n\leq 4$.
\rm

\medskip

\noindent{\bf Proof of Proposition 3.1}. \, By hypothesis, we have
a decomposition as in $(3.2)$, and an $\ell$-adic representation
$\tau'_\ell$ associated to $\eta'$.

As a short digression let us note that if $\eta$ were essentially
self-dual and regular, we could exploit its algebraicity, and by
appealing to [Pic] associate a $3$-dimensional $\ell$-adic
representation to $\eta$. The Proposition $3.1$ will follow in
that case, as in the proof of Proposition $2.1$. However, we
cannot (and do not wish to) assume either that $\eta$ is
essentially self-dual or that it is regular. We have to appeal to
another idea, and here it is.

\medskip

Let $S$ be a finite set of primes including the archimedean and
ramified ones. At any $p$ outside $S$, let $\pi_p$ be represented
by an unordered $(3+r)$-tuple $\{\alpha_1,\dots,\alpha_{3+r}\}$ of
complex numbers, and we may assume that $\eta_p$ (resp. $\eta'_p$)
is represented by $\{\alpha_1,\alpha_2,\alpha_3\}$ (resp.
$\{\alpha_4,\dots,\alpha_{3+r}\}$. It is then straightforward to
check that
$$
L(s,\pi_p; \Lambda^2) \, = \,
L(s,\eta_p^\vee)L(s,\eta_p\times\eta'_p)L(s,\eta';
\Lambda^2).\leqno(3.6)
$$
One can also deduce this as follows. Let $\sigma(\beta)$ denote, for
any irreducible admissible representation $\beta$ of GL$(m,\Q_p)$,
the $m$-dimensional representation of the extended Weil group
$W_{\Q_p}'= W_{\Q_p} \times {\rm SL}(2,\C)$ defined by the local
Langlands correspondence (cf. [HaT], [He]). For any representation
$\sigma$ of $W_{\Q_p}'$ which splits as a direct sum $\tau \oplus
\tau'$, we have
$$
\Lambda^2(\sigma) \, \simeq \, \Lambda^2(\tau)\, \oplus \,
\tau\otimes\tau' \, \oplus \Lambda^2(\tau'),\leqno(3.7)
$$
with $\Lambda^2(\tau')=0$ if $\tau'$ is $1$-dimensional and
$$
\Lambda^2(\tau)\simeq\tau^\vee \quad {\rm when} \quad {\rm
dim}(\tau)=3.\leqno(3.8)
$$
In fact, this shows that the identity (3.6) works at the ramified
primes as well, but we do not need it.

Now, since $\pi^\vee \simeq \eta^\vee \boxplus {\eta'}^\vee$, we
get by putting (3.2) and (3.6) together,
$$
L(s,\pi_p^\vee)L(s,\pi_p\times\eta'_p)L(s,\eta'_p;\Lambda^2) \, =
\,
L(s,\pi_p;\Lambda^2)L(s,{\eta'_p}^\vee)
L(s,\eta'_p\times\eta'_p).\leqno(3.9)
$$
Appealing to Tchebotarev, and using the correspondences
$\pi\leftrightarrow\rho_\ell$ and
$\eta'\leftrightarrow\tau'_\ell$, we obtain the following
isomorphism of $\scG_\Q$-representations:
$$
\rho_\ell^\vee \oplus (\rho_\ell\otimes\tau'_\ell) \oplus
\Lambda^2(\tau'_\ell) \, \simeq \,
\Lambda^2(\rho_\ell)\oplus{\tau'_\ell}^\vee
\oplus\left(\tau'_\ell\otimes\tau'_\ell\right).\leqno(3.10)
$$
Using the decomposition
$$
\tau'_\ell\otimes\tau'_\ell \, \simeq {\rm
sym}^2(\tau'_\ell)\oplus\Lambda^2(\tau'_\ell),
$$
we then obtain $(3.3)$ from $(3.10)$.

\qed

\bigskip

\section{Galois representations attached to regular,
selfdual cusp forms on GL$(4)$}

\medskip

A cusp form $\Pi$ on GL$(m)/F$, $F$ a number field, is said to be
{\it essentially selfdual} iff $\Pi^\vee \simeq \Pi\otimes\lambda$
for an idele class character $\lambda$; it is {\it selfdual} if
$\lambda=1$. We will call such a $\lambda$ a {\it polarization}. Let
us call $\Pi$ {\it almost selfdual} if there is a polarization
$\lambda$ of the form $\mu^2\vert\cdot\vert^t$ for some $t\in\C$ and
a finite order character $\mu$; in this case, one sees that
$(\Pi\otimes\mu)^\vee\simeq \Pi\otimes\mu\vert\cdot\vert^t$, or
equivalently, $\Pi\otimes\mu\vert\cdot\vert^{t/2}$ is selfdual.
Clearly, if $\Pi$ is essentially selfdual, then it becomes, under
base change ([AC]), almost selfdual over a finite cyclic extension
$K$ of $F$.

\medskip

Note that when $\Pi$ is essentially selfdual relative to $\lambda$,
it is immediate that $\lambda_\infty$ occurs in the isobaric sum
decomposition of $\Pi_\infty\boxtimes\Pi_\infty$, or equivalently,
$\sigma(\lambda_\infty)$ is a constituent of
$\sigma(\Pi_\infty)^{\otimes 2}$. This implies that if $\Pi$ is
algebraic, then so is $\Lambda$, and thus corresponds to an
$\ell$-adic character $\lambda_\ell$ of $\scG_\Q$.

Whether or not $\Pi$ is algebraic, we have, for any $S$,
$$
L^S(s,\Pi\times \Pi\otimes\lambda^{-1}) \, = \, L^S(s,\Pi,{\rm
sym}^2\otimes\lambda^{-1})L^S(s,\Pi;\Lambda^2\otimes\lambda^{-1}),\leqno(4.1)
$$
The $L$-function on the left has a pole at $s=1$, since
$\Pi^\vee\simeq\Pi\otimes\lambda$ by hypothesis. Also, neither of
the $L$-functions on the right is zero at $s=1$ ([JS2]).
Consequently, exactly one of the $L$-functions on the right of (4.1)
admits a pole at $s=1$. One says that $\Pi$ is of {\it orthogonal
type} ([Ra3]), resp. {\it symplectic type}, if $L^S(s,\Pi,{\rm
sym}^2\otimes\lambda^{-1})$, resp.
$L^S(s,\Pi,\Lambda^2\otimes\lambda^{-1})$ admits a pole at $s=1$.

\medskip

The following result is a consequence of a synthesis of the results
of a number of mathematicians, and it will be crucial to us in the
next section, while proving Theorem A for $n=5$.

\medskip

\noindent{\bf Theorem B} \, \it Let $\Pi$ be a regular, algebraic
cusp form on GL$(4)/\Q$, which is almost selfdual. Then there exists
a continuous representation
$$
R_\ell: \scG_\Q \, \rightarrow \, {\rm GL}(4,\overline\Q_\ell),
$$
such that
$$
L^S(s,\Pi;\Lambda^2) \, = \, L^S(s,\Lambda^2(R_\ell)),
$$
for a finite set $S$ of primes containing the ramified ones.
Moreover, if $\Pi$ is of orthogonal type, we can show that $R_\ell$
and $\Pi$ are associated, i.e., have the same degree $4$
$L$-functions (outside $S$).

\rm

\medskip

When $\Pi$ admits a discrete series component $\Pi_p$ at some
(finite) prime $p$, a stronger form of this result, and in fact its
generalization to GL$(n)/\Q$, is due to Clozel ([C$\ell$2]). But in
the application considered in the next section, we will not be able
to satisfy such a ramification assumption at a finite place.

In the {\it orthogonal case}, $\Pi$ descends by the work of
Ginzburg-Rallis-Soudry (cf. [So]) to define a regular cusp form
$\beta$ on the split SGO$(4)/\Q$, which is given by a pair
$(\pi_1,\pi_2)$ of regular cusp forms on GL$(2)/\Q$. By Deligne,
there are $2$-dimensional (irreducible) $\ell$-adic representations
$\tau_{1,\ell},\tau_{2,\ell}$, with
$\tau_{j,\ell}\leftrightarrow\pi_j$, $j=1,2$. This leads to the
desired $4$-dimensional $\overline\Q_\ell$-representation $R_\ell:=
\tau_{1,\ell}\otimes\tau_{2,\ell}$ of $\scG_\Q$ associated to $\Pi$,
such that
$$
L^S(s,R_\ell) \, = \, L^S(s,\Pi).\leqno(4.2)
$$
It may be useful to notice that since the polarization is a square
(under the almost selfduality assumption), the associated Galois
representation takes values in SGO$(4,\overline\Q_\ell)$, which is
the connected component of GO$(4,\overline\Q_\ell)$, with quotient
$\{\pm 1\}$. In the general case, not needed for this article,
$R_\ell$ will need to be either of the type above or of {\it Asai
type} (see [Ra4]), associated to a $2$-dimensional
$\overline\Q_\ell$-representation Gal$(\overline\Q/K)$ for a
quadratic extension $K/\Q$.

\medskip

In the (more subtle) {\it symplectic case}, this Theorem is proved
in my joint work [Ra-Sh] with F.~Shahidi. We will start with a
historical comment and then sketch the proof (for the benefit of the
reader). Some years ago, Jacquet, Piatetski-Shapiro and Shalika
announced a theorem, asserting that one could descend any $\Pi$ (of
symplectic type on GL$(4)/\Q$) to a generic cusp form $\beta$ on
GSp$(4)/\Q$ with the same (incomplete) degree $4$ $L$-functions.
Unfortunately, this work was never published, except for a part of
it in [JSh2]. In [Ra-Sh], Shahidi and I provide an alternate,
somewhat more circuitous route, yielding something slightly weaker,
but sufficient for many purposes. Here is the idea. We begin by
considering the twist $\Pi_0:=\Pi\otimes\mu\vert\cdot\vert^{t/2}$
instead of $\Pi$, to make the polarization is trivial, i.e., so that
$\Pi_0$ has parameter in Sp$(4,\C)$. Using the {\it backwards
lifting} results of [GRS] (see also [So]), we get a generic cusp
form $\Pi'$ on the split SO$(5)/\Q$, such that $\Pi_0\to \Pi'$ is
functorial at the archimedean and unramified places. Using the
isomorphism of PSp$(4)/\Q$ with SO$(5)/\Q$, we may lift $\Pi'$ to a
generic cusp form $\tilde{\Pi'}$ on Sp$(4)/\Q$. By a suitable
extension followed by induction, we can associate a generic cusp
form $\Pi_1$ on GSp$(4)/\Q$, such that the following hold:
\newpage
\noindent{$(4.3)$}
\begin{enumerate}
\item[(i)]The archimedean parameter of
$\Pi_2:=\Pi_1\otimes\mu^{-1}\vert\cdot\vert^{-t/2}$ is algebraic and
regular; \, and
\item[(ii)]$L(s,\Pi_{2,p};\Lambda^2) \, = \,
L(s, \Pi_p;\Lambda^2)$ at any prime $p$ where $\Pi$ is unramified.
\end{enumerate}

We in fact deduce a stronger statement in [Ra-Sh], involving also
the ramified primes, but it is not necessary for the application
considered in this paper. To continue, part $(i)$ of $(4.3)$ implies
that $\Pi_2$ contributes to the (intersection) cohomology of (the
Baily-Borel-Satake compactification over $\Q$ of) the
$3$-dimensional Shimura variety $Sh_K/\Q$ associated to GSp$(4)/\Q$,
relative to a compact open subgroup $K$ of GSp$(4,\A_f)$; $Sh_K$
parametrizes principally polarized abelian surfaces with level
$K$-structure. Now by appealing to the deep (independent) works of
G.~Laumon ([Lau1,2]) and R.~Weissauer ([Wei]), one gets a continuous
$4$-dimensional $\ell$-adic representation $R_\ell$ of $\scG_\Q$
such that
$$
L^S(s,\Pi_2) \, = \, L^S(s,R_\ell).\leqno(4.4)
$$
The assertion of Theorem B now follows by combining $(4.3)(ii)$ and
$(4.4)$.

\qed

\bigskip

\section{Two useful Lemmas on cusp forms on GL${\bf(4)}$}

\medskip

Let $F$ be a number field and $\eta$ a cuspidal automorphic
representation of GL$(4,\A_F)$, where $\A_F:=\A \otimes_\Q F$ is the
Adele ring of $F$. Denote by $\omega_\eta$ the central character of
$\eta$.

\medskip

First let us recall (see $(1.2.11)$) that by a difficult {\it
theorem of H.~Kim} ([K]), there is an isobaric automorphic form
$\Lambda^2(\eta)$ on GL$(6)/\Q$ such that
$$
L(s,\Lambda^2(\eta)) \, = \, L(s,\eta; \Lambda^2).\leqno(5.1)
$$

\medskip

\noindent{\bf Lemma 5.2} \, \it $\Lambda^2(\eta)$ is essentially
selfdual. In fact
$$
\Lambda^2(\eta)^\vee \, \simeq \,
\Lambda^2(\eta)\otimes\omega_\eta^{-1}\leqno(5.3)
$$
\rm

\medskip

\noindent{\bf Proof}. \, Thanks to the strong multiplicity one
theorem for isobaric automorphic representations ([JS1]), it
suffices to check this at the primes $p$ where $\eta$ is unramified.
Fix any such $p$, and represent the semisimple conjugacy class
$A_p(\eta)$ by $[a,b,c,d]$. Then it is easy to check that
$$
A_p(\Lambda^2(\eta)) \, = \, \Lambda^2(A_p(\eta)) \, = \,
[ab,ac,ad,bc,bd,cd].\leqno(5.4)
$$
Since for any automorphic representation $\Pi$, the unordered tuple
representing $A_p(\Pi^\vee)$ consists of the inverses of the
elements of tuple representing $A_p(\pi)$, and since
$A_p(\omega_\eta) = [abcd]$, we have
$$
A_p(\Lambda^2(\eta)^\vee\otimes \omega_\eta) \, = \,
[(ab)^{-1},(ac)^{-1},(ad)^{-1},(bc)^{-1},(bd)^{-1},(cd)^{-1}]
\otimes[abcd], \leqno(5.5)
$$
which is none other than $A_p(\Lambda^2(\eta))$. The isomorphism
$(5.3)$ follows.

\qed

\medskip

\noindent{\bf Lemma 5.6} \, \it Let $\eta$ be a cusp form on
GL$(4)/F$ with trivial central character. Suppose $\eta^\vee
\not\simeq \eta$. Then there are infinitely many primes $P$ in
$\scO_F$ where $\eta_P$ is unramified such that $1$ is not an
eigenvalue of the conjugacy class $A_P(\Lambda^2(\eta))$ of
$\Lambda^2(\eta_P)$.\rm

\medskip

\noindent{\bf Proof of Lemma 5.6}. \, Since
$\eta^\vee\not\simeq\eta$, there exist, by the strong multiplicity
one theorem, infinitely many unramified primes $P$ where
$\eta_P^\vee \not\simeq \eta_P$. Pick any such $P$, write
$$
A_P(\eta) \, = \, [a,b,c,d], \quad {\rm with} \quad abcd=1.
$$
The fact that $\eta_P^\vee \not\simeq \eta_P$ implies that the set
$\{a,b,c,d\}$ is not stable under inversion. Hence one of its
elements, which we may assume to be $a$ after renaming, satisfies
the following:
$$
a \, \notin \, \{a^{-1}, b^{-1}, c^{-1}, d^{-1}\}.
$$
Equivalently,
$$
1 \, \notin \{a^2, ab, ac, ad\}.
$$
On the other hand, we have $(5.4)$, using which we conclude that the
only way $1$ can be in this set (attached to $\Lambda^2(\eta_P)$) is
to have either $bc$ or $bd$ or $cd$ to be $1$. But if $bc=1$ (resp.
$bd=1$), since $abcd=1$, we must have $ad=1$ (resp. $ac=1$), which
is impossible. Similarly, if $cd=1$, we are forced to have $ab=1$,
which is also impossible.

\qed

\medskip

\section{Finale}

\medskip

Let $\rho_\ell, \pi$ be as in Theorem A. In view of Corollary $3.5$,
we may assume from henceforth that $n=5$, and that $\pi$ is
algebraic and regular. Suppose $\pi$ is not cuspidal. In view of
Corollary 2.3 and Corollary 3.4, we must then have the decomposition
$$
\pi \, \simeq \, \eta \, \boxplus \nu,\leqno(6.1)
$$
where $\eta$ is an algebraic, regular cusp form on GL$(4)/\Q$ and
$\nu$ an algebraic Hecke character, with associated $\ell$-adic
character $\nu_\ell$.

\medskip

Note that Theorem A needs to be proved under {\it either} of two
hypotheses. To simplify matters a bit, we will make use of the
following:

\medskip

\noindent{\bf Lemma 6.2} \, \it There is a character $\nu_0$ with
$\nu_0^2=1$ such that for $\mu=\nu_0\nu^{-1}$, if $\pi$ is almost
selfdual, then so is $\pi\otimes\mu^{-1}$. \rm

\medskip

\noindent{\bf Proof of Lemma} \, When $\pi$ is almost selfdual,
there exists, by definition, an idele class character $\mu$ such
that $\pi\otimes \mu$ is selfdual. But this implies, thanks to
$(6.1)$ and the cuspidality of $\eta$, that $\eta\otimes\mu$ is
selfdual and $\mu\nu$ is $1$ or quadratic. We are done by taking
$\nu_0=\mu\nu$.

\qed

\medskip

Consequently, we may, and we will, replace $\pi$ by $\pi\otimes\mu$,
$\rho_\ell$ by $\rho_\ell\otimes\mu_\ell$, $\eta$ by
$\eta\otimes\mu$ and $\nu$ by $\nu\mu$, without jeopardizing the
nature of either of the hypotheses of Theorem A. In fact, the {\it
first hypothesis simplifies to assuming that $\pi$ is selfdual}.
Moreover,
$$
\nu^2 \, = \, 1.\leqno(6.3)
$$

\medskip

\noindent{\bf Proof of Theorem A when $\pi$ is almost selfdual}: \,

\medskip

We have to rule out the decomposition $(6.1)$, which gives (for any
finite set $S$ of places containing the ramified and unramified
ones):
$$
L^S(s, \pi) \, = \, L^S(s,\eta)L^S(s,\nu)\leqno(6.4)
$$
As noted above, we may in fact assume that $\pi$ is selfdual and
that $\nu^2=1$. Then the cusp form $\eta$ will also be selfdual and
algebraic. We may then apply Theorem B and conclude the existence of
a $4$-dimensional, semisimple $\ell$-adic representation $\tau_\ell$
associated to $\eta$. Then, expanding $S$ to include $\ell$, we see
that $(6.2)$ implies, in conjunction with the associations
$\pi\leftrightarrow \rho_\ell, \, \eta\leftrightarrow\tau_\ell$,
$$
L^S(s,\rho_\ell) \, = \, L^S(s,\tau_\ell)L^S(s,\nu_\ell).\leqno(6.5)
$$
By Tchebotarev, this gives the isomorphism
$$
\rho_\ell \, \simeq \, \tau_\ell \oplus \nu_\ell,\leqno(6.6)
$$
which contradicts the irreducibility of $\rho_\ell$.

\qed

\medskip

\noindent{\bf Proof of Theorem A for general regular $\pi$}: \,

\medskip

Suppose we have the decomposition $(6.1)$. Again, we may assume that
$\nu^2=1$.

\medskip

Let $\omega=\omega_\pi$ denote the central character of $\pi$. Then
from $(6.1)$ we obtain
$$
\omega \, = \, \omega_\eta\nu.\leqno(6.7)
$$

\medskip

\noindent{\bf Proposition 6.8} \, \it Assume the decomposition
(6.1), and denote by $\omega$ the central character of $\pi$ with
corresponding $\ell$-adic character $\omega_\ell$.
\begin{enumerate}
\item[(a)]We have the identity
$$
L^S(s,\pi;\Lambda^2)L^S(s,\pi^\vee\otimes\omega\nu)\zeta^S(s) \, =
\,
L^S(s,\pi^\vee;\Lambda^2\otimes\omega)L^S(\pi\otimes\nu)L^S(s,\omega\nu).
$$
\item[(b)]There is an isomorphism of $\scG_\Q$-modules
$$
\Lambda^2(\rho_\ell)\oplus(\rho_\ell^\vee\otimes \omega_\ell)\oplus
{\underline 1} \, \simeq \,
\left(\Lambda^2(\rho_\ell^\vee)\otimes\omega_\ell\nu_\ell\right)\oplus
(\rho_\ell\otimes\nu_\ell)\oplus \omega_\ell\nu_\ell.
$$
\end{enumerate}
\rm

\medskip

\noindent{\bf Proof of Proposition 6.8}. \, {\bf (a)} \, It is
immediate, by checking at each unramified prime, that
$$
L^S(s,\pi;\Lambda^2) \, = \,
L^S(s,\Lambda^2(\eta))L^S(\eta\otimes\nu),\leqno(6.9)
$$
and (since $\nu=\nu^{-1}$)
$$
L^S(s,\pi^\vee;\Lambda^2) \, = \,
L^S(s,\Lambda^2(\eta^\vee))L^S(\eta^\vee\otimes\nu).\leqno(6.10)
$$
Since $\omega=\omega_\eta\nu$, we get from Lemma $5.2$ that
$\Lambda^2(\eta^\vee)$ is isomorphic to
$\Lambda^2(\eta)\otimes\omega^{-1}\nu$. Twisting $(6.10)$ by
$\omega\nu$, and using the fact that
$$
L^S(\eta^\vee\otimes\omega)=L^S(s,\pi^\vee\otimes\omega)/L^S(s,\omega\nu),\leqno(6.11)
$$
we obtain
$$
L^S(s,\pi^\vee;\Lambda^2\otimes\omega\nu)L^S(s,\omega\nu) \, = \,
L^S(s,\Lambda^2(\eta))L^S(s,\pi^\vee\otimes\omega).\leqno(6.12)
$$
Similarly, using $(6.9)$ and the fact that
$$
L^S(s,\eta\otimes\nu) \, = \, L^S(s,\pi\otimes\nu)/\zeta^S(s),
$$
we obtain the identity
$$
L^S(s,\pi;\Lambda^2)\zeta^S(s) \, = \,
L^S(s,\Lambda^2(\eta))L^S(s,\pi\otimes\nu).\leqno(6.13)
$$
The assertion of part (a) of the Proposition now follows by
comparing $(6.12)$ and $(6.13)$.

\medskip

{\bf (b)} \, Follows from part (a) by applying Tchebotarev, since
$\rho_\ell\leftrightarrow \pi$.

\qed

\medskip

\noindent{\bf Proposition 6.14} \, We have
\begin{enumerate}
\item[(a)]\, $\omega\nu \, = \, 1.$
\item[(b)] \, $\rho_\ell^\vee \, \simeq \, \rho_\ell$.
\end{enumerate}
\rm

\medskip

\noindent{\bf Proof of Proposition 6.14}. \, ${\bf (a)}$ \,  Since
$\rho_\ell$ is irreducible of dimension $5$, it cannot admit a
one-dimensional summand, and hence part $(b)$ of Proposition $6.8$
implies that either $\omega_\ell\nu_\ell =1$ or
$$
\omega_\ell\nu_\ell \, \subset \, \Lambda^2(\rho_\ell).
$$
Since the first case gives the assertion, let us assume that we are
in the second case. But then, again since $\rho_\ell$ is
irreducible, and since $\Lambda^2(\rho_\ell)$ is a summand of
$\rho_\ell\otimes \rho_\ell$, we must have
$$
\rho_\ell^\vee \, \simeq \, \rho_\ell\otimes
(\omega_\ell\nu_\ell)^{-1}.
$$
In other words, $\rho_\ell$ is essentially selfdual in this case.
Then so is $\pi$. More explicitly, we have (since $\nu^2=1$)
$$
\eta^\vee\boxplus\nu \, \simeq \, \pi^\vee \, \simeq \,
\pi\otimes\omega^{-1}\nu \, \simeq \, \eta\otimes\omega^{-1}\nu
\boxplus \omega^{-1}.
$$
As $\eta$ is cuspidal, this forces the identity
$$
\nu \, = \, \omega^{-1}.
$$
Done.

\medskip

${\bf (b)}$ \, Thanks to part $(a)$ (of this Proposition), we may
rewrite part $(b)$ of Proposition $6.8$ as giving the isomorphism of
$\scG_\Q$-modules
$$
\Lambda^2(\rho_\ell)\oplus(\rho_\ell^\vee\otimes\omega_\ell) \,
\simeq \, \left(\Lambda^2(\rho_\ell^\vee)\right)\oplus
(\rho_\ell\otimes\nu_\ell).\leqno(6.15)
$$
We now need the following:

\medskip

\noindent{\bf Lemma 6.16} \, \it Suppose $\rho_\ell$ is not
selfdual. Then
$$ \rho_\ell\otimes \nu_\ell
\, \not\subset \, \Lambda^2(\rho_\ell).\leqno(6.17)
$$
\rm

\medskip

\noindent{\bf Proof of Lemma $6.16$}. \, The hypothesis on
$\rho_\ell$ implies that $\pi$ is not selfdual, and since
$\pi=\eta\boxplus \nu$, $\eta$ is not selfdual either.

Suppose $(6.17)$ is false.. Then, since $\rho_\ell \leftrightarrow
\pi$, we must have
$$
A_p(\pi\otimes\nu) \, \subset \, A_p(\pi; \Lambda^2), \, \, \forall
\, p\notin S,\leqno(6.18)
$$
for a {\it finite} set $S$ of primes. But we also have
$$
A_p(\pi\otimes\nu) \, = \, A_p(\eta\otimes\nu) \oplus
\underline{1}\leqno(6.19)
$$
and
$$
A_p(\pi; \Lambda^2) \, = \, A_p(\Lambda^2(\eta)) \oplus
A_p(\eta\otimes\nu).\leqno(6.20)
$$
Substituting $(6.19)$ and $(6.20)$ in $(6.18)$, we obtain
$$
\underline{1} \, \subset \, A_p(\Lambda^2(\eta)) \, \, \forall \,
p\notin S.\leqno(6.21)
$$

On the other hand, since $\eta$ is a non-selfdual cusp form on
GL$(4)/\Q$ of trivial central character, we may apply Lemma $5.6$
with $F=\Q$, and conclude that there is an {\it infinite set of
primes} $T$ such that
$$
\underline{1} \, \not\subset \, A_p(\Lambda^2(\eta)) \, \, \forall
\, p\in T,\leqno(6.22)
$$
which contradicts $(6.21)$, proving the Lemma.

\qed

\medskip

In view of the identity $(6.15)$ and Lemma $6.16$, we have now
proved all of Proposition $6.14$.

\qed

\medskip

We are also done with the proof of Theorem A because $\pi$ is
selfdual when the decomposition $(6.1)$ holds, thanks to the
irreducibility of $\rho_\ell$, and the selfdual case has already
been established (using the algebraic regularity of $\pi$).

\qed

\vskip 0.2in

\section*{\bf Bibliography}

\begin{description}

\item[{[AC]}] J.~Arthur and L.~Clozel, \emph{Simple Algebras, Base
Change and the Advanced Theory of the Trace Formula}, Ann. Math.
Studies {\bf 120}, Princeton, NJ (1989).

\item[{[BHR]}] D.~Blasius, M.~Harris and D.~Ramakrishnan,
\emph{Coherent cohomology, limits of discrete series, and Galois
conjugation},  Duke Math. Journal {\bf 73} (1994), no. 3, 647--685.

\item[{[C$\ell$1]}] L.~Clozel, \emph{Motifs et formes automorphes},
in {\it Automorphic Forms, Shimura varieties, and $L$-functions},
vol. I, 77--159, Perspectives in Math. {\bf 10} (1990)

\item[{[C$\ell 2$]}] L.~Clozel, \emph{Repr\'esentations galoisiennes
associ\'ees aux repr\'esentations automorphes autoduales de ${\rm
GL}(n)$}, Publ. Math. IHES {\bf 73}, 97--145 (1991).

\item[{[CoPS1]}] J.~Cogdell and I.~Piatetski-Shapiro,
\emph{Converse theorems for ${\rm GL}\sb n$ II},  J. Reine Angew.
Math. {\bf 507}, 165--188 (1999).

\item[{[CoPS2]}] J.~Cogdell and I.~Piatetski-Shapiro,
\emph{Remarks on Rankin-Selberg convolutions}, in {\sl Contributions
to automorphic forms, geometry, and number theory},  255--278, Johns
Hopkins Univ. Press, Baltimore, MD (2004).

\item[{[De]}] P.~Deligne, \emph{Formes modulaires et Rep\'esentations
$\ell$-adiques}, (1972).

\item[{[De-S]}] P.~Deligne and J.-P.~Serre, \emph{Formes modulaires
de poids $1$}, Ann. Sci. École Norm. Sup. (4) {\bf 7}, 507--530
(1975).

\item[{[FoM]}] J.-M.~Fontaine and B.~Mazur,
\emph{Geometric Galois representations}, in
{\it Elliptic curves,
modular forms, and Fermat's last theorem}, 41--78, Ser. Number
Theory, I, Internat. Press, Cambridge, MA (1995).

\item[{[GJ]}] S.~Gelbart and H.~Jacquet,
\emph{A relation between automorphic
representations of GL$(2)$ and GL$(3)$},
 Ann. Scient. \'Ec. Norm. Sup. (4)
{\bf 11} (1979), 471--542.

\item[{[GRS]}] D.~Ginzburg, S.~Rallis and D.~Soudry,
\emph{ On explicit lifts of cusp forms from ${\rm GL}\sb m$ to
classical groups}, Ann. of Math. (2) {\bf 150}, no. 3, 807--866
(1999).

\item[{[HaT]}] M.~Harris and R.~Taylor,
\emph{On the geometry and cohomology of some simple Shimura
varieties}, with an appendix by V.G.~Berkovich, Annals of Math
Studies {\bf 151}, Princeton University Press, Princeton, NJ (2001)

\item[{[He]}] G.~Henniart, \emph{Une preuve simple des conjectures de Langlands pour
${\rm GL}(n)$ sur un corps $p$-adique}, Invent. Math. {\bf 139}, no. 2, 439--455
(2000).

\item[{[J]}] H.~Jacquet, \emph{Principal $L$-functions of the linear group},
in {\sl Automorphic forms, representations and $L$-functions}, Proc.
Sympos. Pure Math. {\bf 33}, Part 2, 63--86, Amer. Math. Soc.,
Providence, R.I. (1979).

\item[{[JPSS]}] H.~Jacquet, I.~Piatetski-Shapiro and J.A.~Shalika,
\emph{Rankin-Selberg convolutions}, Amer. J of Math. {\bf 105} (1983), 367--464.

\item[{[JS1]}] H.~Jacquet and J.A.~Shalika,
\emph{Euler products and
the classification of automorphic forms} I \& II, Amer. J of
Math. {\bf 103} (1981), 499--558 \& 777--815.

\item[{[JS2]}] H.~Jacquet and J.A.~Shalika, \emph{Exterior square $L$-functions},
in {\it Automorphic forms,
Shimura varieties, and $L$-functions}, Vol. II, 143--226,
Perspectives in Math. {\bf 11} (1990), Academic Press, Boston, MA.

\item[{[K]}] H.~Kim, \emph{Functoriality of the exterior square of GL$_4$ and the
symmetric fourth of GL$_2$}, with Appendix 1 by D.~Ramakrishnan and
Appendix 2 by Kim and P.~Sarnak, J. Amer. Math. Soc. {\bf 16}, no.
1, 139--183 (2003).

\item[{[KSh1]}] H.~Kim and F.~Shahidi,
\emph{Functorial products for GL$(2) \times $GL$(3)$ and the
symmetric cube for GL$(2)$}, With an appendix by Colin J. Bushnell
and Guy Henniart, Annals of Math. (2) {\bf 155}, no. 3, 837--893
(2002).

\item[{[KSh2]}] H.~Kim and F.~Shahidi, \emph{Cuspidality of
symmetric powers with applications}, Duke Math. J. {\bf 112},  no.
1, 177--197 (2002).

\item[{[La1]}] R.P.~Langlands, \emph{On the classification of
irreducible representations of real algebraic groups}, in {\it
Representation theory and harmonic analysis on semisimple Lie
groups}, 101--170, Math. Surveys Monographs {\bf 31}, AMS,
Providence, RI (1989).

\item[{[La2]}] R.P.~Langlands, \emph{Automorphic representations,
Shimura varieties, and motives. Ein M\"archen}, Proc. symp. Pure
Math {\bf 33}, ed. by A. Borel and W. Casselman, part 2, 205--246,
Amer. Math. Soc., Providence (1979).

\item[{[La3]}] R.P.~Langlands, \emph{On the notion of an automorphic
representation}, Proc. symp. Pure Math {\bf 33}, ed. by A. Borel and
W. Casselman, part 2, 189--217, Amer. Math. Soc., Providence (1979).

\item[{[Lau1]}] G.~Laumon,  \emph{Sur la cohomologie \`a supports
compacts des vari\'et\'es de Shimura
pour ${\rm GSp}(4)/\Q$}, Compositio Math. {\bf 105} (1997), no. 3,
267--359.

\item[{[Lau2]}] G.~Laumon, \emph{Fonctions z\'etas des
vari\'et\'es de Siegel de dimension trois},
in {\sl Formes automorphes II: le cas du groupe GSp(4)}, Edited by
J. Tilouine, H. Carayol,  M. Harris, M.-F. Vigneras, Asterisque {\bf
302}, Soc. Math. France {\it Ast\'erisuqe} (2006).

\item[{[MW]}] C.~Moeglin and J.-L.~Waldspurger, \emph{Poles des
fonctions $L$ de paires pour GL$(N)$},
Appendice, Ann. Sci. \'Ecole Norm. Sup. (4) {\bf 22} (1989),
667-674.

\item[{[Pic]}] \emph{Zeta Functions of Picard Modular Surfaces},
edited by R.P.~Langlands and D.~Ramakrishnan, CRM Publications,
Montr\'eal (1992).

\item[{[Ra1]}] D.~Ramakrishnan, \emph{Pure motives and automorphic
forms}, in {\it Motives},
(1994) Proc. Sympos. Pure Math. 55, Part 2, AMS, Providence, RI,
411--446.

\item[{[Ra2]}] D.~Ramakrishnan, \emph{Modularity of the
Rankin-Selberg $L$-series, and
Multiplicity one for SL$(2)$}, Annals of Mathematics
{\bf 152} (2000), 45--111.

\item[{[Ra3]}] D.~Ramakrishnan, \emph{Modularity of solvable Artin
representations of GO$(4)$-type}, IMRN {\bf 2002}, No. {\bf 1}
(2002), 1--54.

\item[{[Ra4]}] D.~Ramakrishnan, \emph{Algebraic cycles on
Hilbert modular fourfolds
and poles of $L$-functions}, in {\it Algebraic groups and
arithmetic}, 221--274, Tata Inst. Fund. Res., Mumbai (2004).

\item[{[Ra5]}] D.~Ramakrishnan, \emph{Irreducibility of $\ell$-adic
associated to regular cusp forms on GL$(4)/\Q$}, preprint (2004),
being revised.

\item[{[Ra-Sh]}] D.~Ramakrishnan and F.~Shahidi, \emph{Siegel
modular forms of genus $2$
attached to elliptic curves}, preprint, submitted (2006)

\item[{[Ra-W]}] D.~Ramakrishnan and S.~Wang, \emph{A cuspidality
criterion for the functorial product on $\rm GL(2)\times GL(3)$
with
a cohomological application}, IMRN {\bf 2004}, No. {\bf 27},
1355--1394.

\item[{[Ri]}]K.~Ribet, \emph{Galois representations attached to
eigenforms with Nebentypus}, in {\sl Modular functions of one
variable V}, pp. 17--51, Lecture Notes in Math. {\bf 601}, Springer,
Berlin (1977).

\item[{[Se]}] J.-P.~Serre, \emph{Abelian $\ell$-adic
representations}, Research Notes in Mathematics {\bf 7}, A.K.~Peters
Ltd., Wellesley, MA (1998).

\item[{[Sh1]}] F.~Shahidi, \emph{On the Ramanujan conjecture and the
finiteness of poles for certain $L$-functions}, Ann. of Math. (2) {\bf 127}
(1988), 547--584.

\item[{[Sh2]}] F.~Shahidi, \emph{A proof of the Langlands conjecture on
Plancherel measures; Complementary series for $p$-adic groups}, Ann. of
Math. {\bf 132} (1990), 273-330.

\item[{[So]}]D.~Soudry,  \emph{On Langlands functoriality from classical groups
to ${\rm GL}\sb n$}, in {\it Automorphic forms I}, Ast\'erisque {\bf
298}, 335--390 (2005).

\item[{[T]}] J.~Tate, \emph{Les conjectures de Stark sur les
fonctions $L$ d'Artin en $s=0$}, Lecture notes edited by D.~Bernardi
and N.~Schappacher, Progress in Mathematics {\bf 47} (1984),
Birkh\"auser, Boston, MA.

\item[{[Wei]}] R.~Weissauer, \emph{Four dimensional Galois
representations}, in {\sl Formes automorphes II: le cas du groupe
GSp(4)}, Edited by J. Tilouine, H. Carayol,  M. Harris, M.-F.
Vigneras, Asterisque {\bf 302}, Soc. Math. France {\it Ast\'erisuqe}
(2006).
\bigskip

\end{description}

\bigskip

Dinakar Ramakrishnan

253-37 Caltech

Pasadena, CA 91125, USA.

dinakar@caltech.edu

\bigskip

\vskip 0.2in

\end{document}